\documentclass[11pt]{amsart}

\usepackage[T1]{fontenc}
\usepackage{lmodern}
\usepackage{microtype}
\usepackage{amsmath}
\usepackage{amsthm}
\usepackage{amsfonts}
\usepackage{amssymb}
\usepackage{mathtools}
\usepackage{esint}
\usepackage{geometry}
\usepackage{needspace}
\usepackage{hyperref}
\usepackage[nameinlink,capitalise,noabbrev]{cleveref}
\usepackage{mathrsfs}

\hypersetup{
	hidelinks,
	pdftitle={Local Aronson--B\'enilan estimates for the porous medium and fast diffusion equations on collapsing manifolds},
	pdfauthor={Jinpeng Miao and Bin Shen}
}

\newcommand{\vol}{\operatorname{vol}}
\newcommand{\Ric}{\operatorname{Ric}}
\newcommand{\supp}{\operatorname{supp}}
\numberwithin{equation}{section}

\theoremstyle{plain}
\newtheorem{theorem}{Theorem}[section]
\newtheorem{proposition}[theorem]{Proposition}
\newtheorem{lemma}[theorem]{Lemma}
\newtheorem{corollary}[theorem]{Corollary}

\theoremstyle{definition}

\theoremstyle{remark}

\title[Aronson--B\'enilan estimates for PME and FDE]
{Local Aronson--B\'enilan estimates for the porous medium and fast diffusion equations on collapsing manifolds}

\author{Jinpeng Miao}
\address{
	School of Mathematics, Southeast University,
	Nanjing 211189, Jiangsu, China}
\email{jmiao@seu.edu.cn}

\author{Bin Shen*}
\address{
	School of Mathematics, Southeast University,
	Nanjing 211189, Jiangsu, China}
\email{shenbin@seu.edu.cn}
\thanks{*Corresponding author. The second author is supported partially by the Jiangsu Provincial Scientific Research Center of Applied Mathematics (No. BK20233002), and partially by the NNSFC (Nos. 12271093, 12301032).}


\begin{document}
	% ======================================================================
	\begin{abstract}
		We establish local Aronson--B\'enilan estimates for positive solutions to the porous medium and fast diffusion equations on complete noncompact Riemannian manifolds with a small scale-invariant integral bound on the negative part of the Ricci curvature in $L^p$, where $p>n/2$. No uniform lower bound on the volume of unit balls is assumed, and hence the results apply to manifolds that may collapse. For the fast diffusion equation, we obtain the estimate in the range $1-2/n<m<1$.		
		The auxiliary function method, which was designed for the heat equation under integral curvature conditions, is not directly applicable to those two nonlinear equations. The pressure operator in both cases contains a solution-dependent coefficient, which produces a residual term involving $\partial_t J$ in the evolution inequality. Its treatment requires a positive lower bound for $J$, and moreover, control of $-\partial_t\log J$. We obtain these bounds by linearizing the auxiliary equation, applying Duhamel's formula, and combining local Dirichlet heat kernel estimates with a singular Volterra inequality. Additionally, in the fast diffusion case, the negativity of $v$ changes the sign structure of the evolution inequality, so a separate coercive estimate is required to close the argument. Local Aronson--B\'enilan estimates follow for both diffusion regimes and integrating them along space-time curves gives the corresponding Harnack inequalities.
	\end{abstract}
	
	
	\maketitle

	% ======================================================================
	\section{Introduction}

	It is well known that on a Riemannian manifold with a pointwise lower bound on the Ricci curvature, positive solutions of parabolic equations satisfy sharp Li--Yau gradient estimates \cite{LiYau}. Such estimates have proven to be a powerful tool in geometric analysis, with applications ranging from Harnack inequalities to entropy formulae. However, the pointwise curvature condition $\operatorname{Ric} \ge -Kg$ is somehow restrictive. It fails on manifolds that collapse, or on those with isolated curvature blow-ups.
	
	For the heat equation, Zhang and Zhu overcame this limitation by
	introducing an integral Ricci curvature condition
	\cite{ZhangZhu2017,MR3802491}. Let $\Ric^-$ denote the negative
	spectral part of the Ricci curvature tensor, and let
	$|\Ric^-|_{\mathrm{HS}}$ denote its Hilbert--Schmidt norm. For
	$p>n/2$, they assumed
	\[
	\sup_{x\in M}
	\left(
	\fint_{B(x,1)}|\Ric^-|_{\mathrm{HS}}^p\,dV
	\right)^{1/p}
	\leq\kappa,
	\]
	where $\kappa=\kappa(n,p)>0$ is sufficiently small. The relation
	between this condition and our convention is explained at the
	beginning of \cref{sec:preliminaries}.
	Under this condition, they constructed an auxiliary function $J$ satisfying
	\begin{equation}\label{1.2}
		J_t - \Delta J + 2WJ + a \frac{|\nabla J|^2}{J} = 0,
	\end{equation}
	where $W$ is a smooth majorant of the negative Ricci curvature magnitude and $a > 1$ is a constant. By carefully analyzing the upper and lower bounds of $J$ under Dirichlet boundary conditions, they obtained a Li--Yau gradient estimate for positive solutions of the heat equation. This technique has since been extended to various nonlinear heat-type equations \cite{Wang2020,WangWang2024,MR4665538}.
	
	Nevertheless, when one attempts to apply the Zhang--Zhu method to the porous medium equation (PME) and the fast diffusion equation (FDE)
	\begin{equation}\label{eq:intro-diffusion}
		u_t=\Delta u^m,\qquad
		\begin{cases}
			m>1 & \text{for  PME},\\
			0<m<1 & \text{for FDE},
		\end{cases}
	\end{equation}
	on collapsing manifolds, a substantial obstacle arises.

	In Euclidean space, Aronson and B\'enilan established in \cite{AB79} the classical estimate for the pressure variable
	\begin{equation}\label{Pressvar}
		v=\frac{m}{m-1}u^{m-1}.
	\end{equation}
	They obtained
	\begin{equation*}
		\Delta v\geq-\frac{C}{t}.
	\end{equation*}
	For the PME, this is equivalent to
	\begin{equation*}
		\frac{|\nabla v|^2}{v}-\frac{v_t}{v}
		\leq\frac{C}{t}.
	\end{equation*}
	On a Riemannian manifold, the Bochner formula introduces the curvature term $\operatorname{Ric}(\nabla v,\nabla v)$. Under a pointwise lower bound $\operatorname{Ric}\geq-Kg$, Lu, Ni, V\'azquez, and Villani adapted the maximum principle argument to obtain local Aronson--B\'enilan estimates for both the PME and the FDE \cite{LNVV09}. More recently, Huang and Shen developed refined gradient estimates via Moser iteration and alternative maximum principle techniques, still under pointwise curvature assumptions \cite{HuangShen2025,HuangShen2026}. For the Witten Laplacian on weighted manifolds, Huang and Li further established gradient estimates and entropy formulae for the porous medium and fast diffusion equations \cite{HuangLi2014}. These results demonstrate that the Aronson--B\'enilan estimate is meaningful and available on manifolds, provided the curvature is sufficiently controlled. In order to get such estimates on collapsing manifolds, the way we undertake in this manuscript is to remove the pointwise curvature assumption and replace it with the more flexible integral curvature condition.
	
	The main difficulty comes from the pressure operator
	\begin{equation}\label{operatorL}
		\mathcal{L}
		=
		\partial_t-(m-1)v\Delta,
	\end{equation}
	which contains the $v$-dependent coefficient $(m-1)v$ in front of the Laplacian. It satisfies
	\begin{equation}\label{eq:v-equation-intro}
		v_t=(m-1)v\Delta v+|\nabla v|^2
	\end{equation}
	on the positivity set of $u$. Consequently,
	\begin{equation*}
		\frac{|\nabla v|^2}{v}-\frac{v_t}{v}
		=-(m-1)\Delta v.
	\end{equation*}
	Thus an upper bound for the expression on the left is equivalent to
	a lower bound (PME) or upper bound (FDE)  for the Laplacian of the pressure. When the auxiliary equation \eqref{1.2} for $J$ is inserted into the evolution formula for the quantity
	\begin{equation*}
		F_{\mathrm P} = \alpha J \frac{|\nabla v|^2}{v} - \frac{v_t}{v},
	\end{equation*}
	a residual term of the form
	\begin{equation*}
		\alpha J_t(1-(m-1)v)\frac{|\nabla v|^2}{v}
	\end{equation*}
	appears. Controlling this term requires not only a lower bound for $J$, but also a control on its logarithmic time derivative $-\partial_t\log J$. In the heat equation setting, where the operator is simply $\partial_t - \Delta$, such a term is absent. For the PME and the FDE, however, this is the central obstruction.
	
	To overcome this obstruction, we develop a refined analysis of the auxiliary function $J$. We linearize \eqref{1.2} via the substitution $J = \zeta^{-1/\beta}$, apply Duhamel's formula, and use local heat kernel estimates that are uniform with respect to the intermediate domain, together with a singular Volterra inequality. This yields both a lower bound for $J$ and a two-sided bound for $\partial_t\log J$. These estimates are precisely what is needed to close the maximum principle argument for the PME and the FDE under integral Ricci curvature bounds.
	
	In contrast to the PME case, the FDE requires an additional argument because its pressure variable satisfies \(v<0\), whereas \(v>0\) for the PME. This sign difference changes the structure of the evolution inequality. We therefore retain the Hessian term and establish a positive-definite quadratic estimate to close the maximum principle argument in the range \(1-2/n<m<1\).

	The main result of this manuscript is a local Aronson--B\'enilan estimate for positive solutions of the PME and the FDE on complete noncompact Riemannian manifolds under the normalized integral curvature assumption \eqref{eq:curvature-assumption}, without requiring any uniform lower bound on the volume of unit balls. Consequently, our results apply to manifolds that may collapse. By integrating the Aronson--B\'enilan estimates along space-time curves, we also obtain integrated and distance-form Harnack inequalities for both equations.
	
	The manuscript is organized as follows. The assumptions, notation, main estimates, and Harnack consequences are presented in \cref{sec:main-results}, which concludes with a comparison with previous results in \cref{subsec:comparison}. \cref{sec:preliminaries} collects the geometric preliminaries, including local volume comparison, Sobolev inequalities, Dirichlet heat kernel bounds, and the singular Volterra inequality. The construction of the auxiliary function $J$ and the proof of the key bounds on $J$ and $-\partial_t\log J$ are given in \cref{sec:auxiliary}. The PME evolution inequalities and the localized maximum principle argument are presented in \cref{sec:pme-evolution}. The corresponding FDE evolution inequalities and localization argument are given in \cref{sec:fde-evolution}. Finally, the integrated Harnack inequalities for both equations are proved in \cref{sec:unified-harnack}.
	\section{Main results}
	\label{sec:main-results}
	
	\subsection{Curvature assumption}
	\label{subsec:assumption}
	
	We work on a complete noncompact Riemannian manifold $(M^n,g)$ with
	$n\geq3$. Fix $p>n/2$ and a single
	sufficiently small constant $\kappa=\kappa(n,p)>0$ for the curvature
	assumptions below.
	
	Following the convention of Dai--Wei--Zhang \cite{DWZ}, we measure
	the negative Ricci curvature by the negative part of its smallest
	eigenvalue. More precisely, let $\lambda_{\min}(\Ric_x)$ denote
	the smallest eigenvalue of the Ricci endomorphism at $x$, and define
	\begin{equation}
		\label{eq:rho-minus}
		\rho_-(x)=\max\{0,-\lambda_{\min}(\Ric_x)\}.
	\end{equation}
	The Hilbert--Schmidt norm of the negative spectral part satisfies
	\begin{equation}
		\label{eq:ricci-norm-comparison}
		\rho_-\leq|\Ric^-|_{\mathrm{HS}}\leq\sqrt n\,\rho_-.
	\end{equation}
	Thus the two curvature quantities are equivalent up to a dimensional
	factor. We use the Dai--Wei--Zhang convention $\rho_-$ throughout.
	Its relation to the Zhang--Zhu condition is explained at the beginning
	of \cref{sec:preliminaries}, and smooth majorants of $\rho_-$ are
	constructed in \cref{subsec:negative-ricci}.

	For $x\in M$ and $r>0$, define
	\begin{equation}
		\label{eq:curvature-scale-invariant}
		k(x,p,r)
		=
		r^2
		\left(
		\fint_{B(x,r)}
		\rho_-^{\,p}\,dV
		\right)^{1/p},
		\qquad
		k(p,r)=\sup_{x\in M}k(x,p,r).
	\end{equation}
	Here $\fint_B=\vol(B)^{-1}\int_B$.
	Throughout the main results, we assume that
	\begin{equation}
		\label{eq:curvature-assumption}
		k(p,1)\leq\kappa.
	\end{equation}

	\subsection{Notation for the main estimate}
	\label{subsec:notation-assumptions}
	\label{subsec:notation}
	We use the exponents
	\begin{equation}
		\label{eq:common-exponents}
		\gamma=\frac{n}{2p},
		\qquad
		\sigma=\frac{n}{2p-n}.
	\end{equation}
	
	In the following definitions, $C=C(n,p)>0$ is chosen so that
	 Corollary~\ref{lem:cutoff}, Proposition~\ref{prop:volterra},
	and Lemma~\ref{lem:kernel-potential} hold simultaneously.
	For $\beta>0$, $0<r\leq1$, and $t>0$, define
	\begin{equation}
\label{eq:Lambda-general}
\Lambda(\beta)=C\beta^\sigma\kappa^{1+\sigma},\quad
\underline J_r^{(\beta)}(t)=2^{-1/\beta}e^{-\Lambda(\beta)t/r^2},\quad
G_r^{(\beta)}(t)=C\kappa r^{-2}(t/r^2)^{-\gamma}e^{\beta\Lambda(\beta)t/r^2}.
\end{equation}

	\subsubsection{Coefficients for the PME}\leavevmode\par\smallskip
	\label{subsubsec:pme-coefficients}
	
For the PME, fix $m>1$ and $0<\alpha<1$, and set
\begin{equation}
\label{eq:c-alpha}
c_\alpha=\frac{\alpha(1-\alpha)}2,\qquad
 a_\alpha^{\mathrm{PME}}=\frac{4+2/(m-1)}{1-\alpha},\qquad
 \beta^{\mathrm{PME}}=a_\alpha^{\mathrm{PME}}-1.
\end{equation}
The PME weights are
\begin{equation}\label{eq:underlineJ-PME}
\underline J_{r,\mathrm P}(t)=\underline J_r^{(\beta^{\mathrm{PME}})}(t),\qquad
G_{r,\mathrm P}(t)=G_r^{(\beta^{\mathrm{PME}})}(t).
\end{equation}
	The first is a lower bound for $J$, and the second controls
	$-\partial_t\log J$.

	\subsubsection{Coefficients for the FDE}\leavevmode\par\smallskip
	\label{subsubsec:fde-coefficients}
	
	For the FDE, assume $1-2/n<m<1$ and fix any $\alpha>1$.
    Set $\mu=1-m$. Choose $\varepsilon$ first to determine $\lambda_m$,
    and then choose $\theta$ so that
	\begin{equation}
		\label{eq:lambda-m-choice}
		\begin{aligned}
			0<\varepsilon&<1-\frac{n\mu}{2},\qquad
			\lambda_m=\frac{2(1-\varepsilon)}{n\mu}>1,\\
			0<\theta&<\min\left\{\frac12,
			\frac{(\lambda_m-1)(\alpha-1)}{\lambda_m\alpha-1}\right\}.
		\end{aligned}
	\end{equation}
Define
\begin{equation}\label{eq:a-FDE}
a^{\mathrm{FDE}}=\frac2\varepsilon+\frac1{\mu\theta},\qquad
\beta^{\mathrm{FDE}}=a^{\mathrm{FDE}}-1.
\end{equation}
The FDE weights are
\begin{equation}\label{eq:underlineJ-FDE}
\underline J_{r,\mathrm F}(t)=\underline J_r^{(\beta^{\mathrm{FDE}})}(t),\qquad
G_{r,\mathrm F}(t)=G_r^{(\beta^{\mathrm{FDE}})}(t).
\end{equation}
	By Lemma~\ref{lem:matrix-M-positive}, these choices make
    $M(A)$ positive definite for every $A\geq\alpha$. For fixed $T>0$
	and $0<r\leq1$, define
\begin{equation}
    \label{eq:delta-effective-theorem}
    \delta=\min\left\{1,
    \min_{\alpha\leq A\leq\alpha/\underline J_{r,\mathrm F}(T)}
    \lambda_{\min}(M(A))\right\}.
\end{equation}

	\subsection{Main estimates and Harnack inequalities}
	\label{subsec:main-estimates}
	The following estimates use the parameters and functions fixed in
	\cref{subsec:notation}.
	
	\subsubsection{The porous medium equation}\leavevmode\par\smallskip
	
	\begin{theorem}[Local Aronson--B\'enilan estimate for the PME]
		\label{thm:local-AB}
		Let $(M^n,g)$ be complete and noncompact with $n\geq3$, and let
		$p>n/2$. Assume
		\[
		k(p,1)\leq\kappa.
		\]
		Let $u \in C(M \times [0,T]) \cap C^\infty(M \times (0,T])$ satisfy
		$u>0$ on $M\times[0,T]$ and solve \eqref{eq:intro-diffusion}
        with $m>1$ on $M\times(0,T]$. Let $v$ be the pressure defined
        in \eqref{Pressvar}. Fix $O\in M$, $0<r\leq1$, and
		$0<\alpha<1$. For $0<t\leq T$, set
		\[
		V_r^{\mathrm{P}}(t)
		=
		\sup_{\overline{B(O,r)}\times(0,t]}v.
		\]
		Then on $B(O, r/2) \times (0,T]$,
		\begin{equation}
			\label{eq:main-estimate-PME}
			\begin{aligned}
				\alpha \underline J_{r,\mathrm{P}}(t) \frac{|\nabla v|^2}{v} - \frac{v_t}{v}
				&\le \frac{n(m-1)}{n(m-1)+\alpha(1+\alpha)\underline J_{r,\mathrm P}(t)}\frac1t\\
				&\quad+\frac{Cm^2(m-1)V_r^{\mathrm P}(t)}{c_\alpha\underline J_{r,\mathrm P}(t)r^2}
                +\frac{C(m-1)G_{r,\mathrm P}(t)V_r^{\mathrm P}(t)}{\sqrt{c_\alpha}}.
			\end{aligned}
		\end{equation}
	\end{theorem}

    Integrating this estimate along space-time curves gives a Harnack inequality.
	Let $x_1,x_2\in B(O,r/2)$ and $0<t_1<t_2\le T$. Assume that a piecewise
	$C^1$ map $\eta$ from $[0,1]$ into $B(O,r/2)$ satisfies
	$\eta(0)=x_1$ and $\eta(1)=x_2$. Set $\xi(s)=(1-s)t_1+st_2$ and
	\begin{equation*}
		\begin{aligned}
			\Psi_r(t)
			={}&\frac{Cm^2(m-1)V_r^{\mathrm P}(t_2)}{c_\alpha\underline J_{r,\mathrm P}(t)r^2}
			+\frac{C(m-1)G_{r,\mathrm P}(t)V_r^{\mathrm{P}}(t_2)}{\sqrt{c_\alpha}}.
		\end{aligned}
	\end{equation*}
	
	\begin{theorem}[Integrated Harnack inequality for the PME]
		\label{thm:integrated-Harnack}
		Under the hypotheses of \cref{thm:local-AB},
		\begin{equation}
			\label{eq:harnack-PME}
			\begin{aligned}
				v(x_1,t_1)
				\leq{}&v(x_2,t_2)\exp\Biggl\{
				\int_{t_1}^{t_2}\left[
				\frac{n(m-1)}{n(m-1)+\alpha(1+\alpha)\underline J_{r,\mathrm P}(t)}\frac1t
				+\Psi_r(t)\right]dt\\
				&\qquad+\frac{1}{4\alpha(t_2-t_1)}
				\int_0^1\frac{|\dot\eta(s)|^2}
				{\underline J_{r,\mathrm{P}}(\xi(s))v(\eta(s),\xi(s))}\,ds
				\Biggr\}.
			\end{aligned}
		\end{equation}
	\end{theorem}

	\begin{corollary}[Distance form of Harnack inequality for the PME]
		\label{cor:distance-harnack}
		Under the hypotheses and notation of
		Theorem~\ref{thm:integrated-Harnack}, let $0 < \rho < r/2$ and define
		\[
		v_{*,\rho} = \min_{\overline{B(O,\rho)} \times [t_1, t_2]} v.
		\]
		Assume that $x_1$ and $x_2$ can be joined inside $B(O,\rho)$ by piecewise $C^1$ curves, and let $d_\rho$ denote the corresponding intrinsic length metric. For $x_1, x_2 \in B(O, \rho)$,
		\begin{equation}
			\label{eq:distance-harnack}
			\begin{aligned}
				v(x_1,t_1)
				\leq{}&v(x_2,t_2)\exp\Biggl\{
				\int_{t_1}^{t_2}\left[
				\frac{n(m-1)}{n(m-1)+\alpha(1+\alpha)\underline J_{r,\mathrm P}(t)}\frac1t
				+\Psi_r(t)\right]dt\\
				&\qquad+\frac{d_\rho(x_1,x_2)^2}
				{4\alpha\underline J_{r,\mathrm{P}}(t_2)v_{*,\rho}(t_2-t_1)}
				\Biggr\}.
			\end{aligned}
		\end{equation}
	\end{corollary}
	
	\subsubsection{The fast diffusion equation}\leavevmode\par\smallskip
	
	For the fast diffusion equation,
	\begin{equation}
		\label{eq:FDE}
		u_t = \Delta u^m, \qquad 0 < m < 1.
	\end{equation}
	
	We work in the supercritical range
	\begin{equation}
		\label{eq:FDE-range-theorem}
		1 - \frac{2}{n} < m < 1.
	\end{equation}
	\begin{theorem}[Local Aronson--B\'enilan estimate for the FDE]
		\label{thm:fde-local-AB}
		Let $(M^n,g)$ be complete and noncompact with $n\geq3$, and let
		$p>n/2$. Assume \eqref{eq:FDE-range-theorem} and
		\[
		k(p,1)\leq\kappa.
		\]
		
		Let $u \in C(M \times [0,T]) \cap C^\infty(M \times (0,T])$ satisfy
		$u>0$ on $M\times[0,T]$ and solve \eqref{eq:FDE} on
		$M\times(0,T]$. Define
		$v=\frac{m}{m-1}u^{m-1}<0$. Fix $O\in M$, $0<r\leq1$, and $\alpha>1$.
		For $0<t\leq T$, set
		\[
		V_r^{\mathrm{F}}(t)
		=
		\sup_{\overline{B(O,r)}\times(0,t]}(-v).
		\]
		
		Choose $\varepsilon$ and $\theta$ as in
		\eqref{eq:lambda-m-choice}, and let $\delta$ be given by
		\eqref{eq:delta-effective-theorem}.
				Then, on $B(O,r/2)\times(0,T]$,
		\begin{equation}
			\label{eq:fde-theorem-estimate}
			\begin{aligned}
				\frac{v_t}{v} - \alpha\frac{\underline J_{r,\mathrm F}(t)}{\underline J_{r,\mathrm F}(T)}\frac{|\nabla v|^2}{v}
				&\le \frac{C}{\delta t} + \frac{Cm^2 V_r^{\mathrm{F}}(t)}{\delta^2 r^2} 
				+ \frac{C(1-m) V_r^{\mathrm{F}}(t)}{\delta r^2} + \frac{C\alpha(1-m)G_{r,\mathrm F}(t)V_r^{\mathrm{F}}(t)}{\delta\underline J_{r,\mathrm F}(T)}.
			\end{aligned}
		\end{equation}
	\end{theorem}
	
	The corresponding Harnack inequality follows by integration. Let
	$x_1,x_2\in B(O,r/2)$ and let $0<t_1<t_2\leq T$. Assume that a
	piecewise $C^1$ map $\eta$ from $[0,1]$ into $B(O,r/2)$ joins $x_1$ to
	$x_2$. Set
	\[
	\xi(s)=(1-s)t_1+st_2.
	\]
	Define
	\begin{equation*}
		\begin{aligned}
			\Phi_r(t)
			={}&\frac{C}{\delta t}
			+\frac{Cm^2V_r^{\mathrm{F}}(t_2)}{\delta^2r^2}
			+\frac{C(1-m)V_r^{\mathrm{F}}(t_2)}{\delta r^2}
			+\frac{C\alpha(1-m)G_{r,\mathrm F}(t)V_r^{\mathrm{F}}(t_2)}{\delta\underline J_{r,\mathrm F}(T)}.
		\end{aligned}
	\end{equation*}
	
	\begin{theorem}[Integrated Harnack inequality for the FDE]
		\label{thm:fde-harnack}
		Under the hypotheses of \cref{thm:fde-local-AB},
		\begin{equation}
			\label{eq:fde-harnack}
			\begin{aligned}
				-v(x_2,t_2)
				\leq{}&(-v(x_1,t_1))\exp\Biggl\{
				\int_{t_1}^{t_2}\Phi_r(t)\,dt
				+\frac{\underline J_{r,\mathrm F}(T)}{4\alpha(t_2-t_1)}
				\int_0^1\frac{|\dot\eta(s)|^2}
				{\underline J_{r,\mathrm{F}}(\xi(s))(-v(\eta(s),\xi(s)))}\,ds
				\Biggr\}.
			\end{aligned}
		\end{equation}
	\end{theorem}

	\begin{corollary}[Distance form of Harnack inequality for the FDE]
		\label{cor:fde-distance-harnack}
		Under the hypotheses and notation of Theorem~\ref{thm:fde-harnack}, let
		$0<\rho<r/2$ and define
		\begin{equation*}
			v_{*,\rho}^{\mathrm{F}}
			=\min_{\overline{B(O,\rho)}\times[t_1,t_2]}(-v).
		\end{equation*}
		Assume that $x_1$ and $x_2$ can be joined inside $B(O,\rho)$ by
		piecewise $C^1$ curves, and let $d_\rho$ denote the intrinsic length distance induced by $g$ inside $B(O,\rho)$. For $x_1,x_2\in B(O,\rho)$,
		\begin{equation}
			\label{eq:fde-distance-harnack}
			\begin{aligned}
				-v(x_2,t_2)
				\leq{}&(-v(x_1,t_1))\exp\Biggl\{
				\int_{t_1}^{t_2}\Phi_r(t)\,dt
				+\frac{\underline J_{r,\mathrm F}(T)d_\rho(x_1,x_2)^2}
				{4\alpha\underline J_{r,\mathrm{F}}(t_2)
					v_{*,\rho}^{\mathrm{F}}(t_2-t_1)}
				\Biggr\}.
			\end{aligned}
		\end{equation}
	\end{corollary}

	\subsection{Comparison with previous results}
\label{subsec:comparison}
\leavevmode\par

\noindent The comparison is organized into three parts, covering bounds for the
auxiliary function, gradient estimates under a pointwise Ricci bound,
and sharpness and coefficient choices.

\smallskip\noindent\textbf{Bounds for the auxiliary function.}\par
\noindent With a common sufficiently large $C=C(n,p)$, the exponent of $e$
in the lower bound of Proposition~\ref{prop:J} is larger than that in
\cite[(2.7)]{ZhangZhu2017},
\[
-C\kappa(\beta\kappa)^\sigma\frac{t}{r^2}
>-C\kappa[1+(\beta\kappa)^\sigma]\frac{t}{r^2}.
\]
This gives a larger lower bound for $J$.
In \cite[(2.12)--(2.13)]{ZhangZhu2017},
volume comparison gives the bound
$Ck(p,r)r^{-2}\leq C\kappa r^{-2}$ for $t-s\geq r^2$.
Our estimate \eqref{eq:potential-estimate} retains the factor
$((t-s)/r^2)^{-\gamma}$ on this interval. This avoids the
constant large time contribution that produces the extra $1$
in the exponent of $e$ in their lower bound.

\smallskip\noindent\textbf{Gradient estimates under a pointwise Ricci bound.}
We replace the integral curvature assumption by a pointwise Ricci lower bound
and compare the resulting estimates with those in
\cite{LNVV09,HuangHuangLi2013,HuangShen2025}.

Let $(M^n,g)$ be complete and noncompact, with $n\geq3$ and
$\Ric\geq-Kg$, $K\geq0$. Let $T>0$ and let $u>0$ solve $u_t=\Delta u^m$
on $M\times(0,T]$, with $u\in C(M\times[0,T])\cap C^\infty(M\times(0,T])$.
Use the pressure \eqref{Pressvar} and set $a_m=n(m-1)/(n(m-1)+2)$.

For the PME, fix $m>1$ and $0<\alpha<1$, and set
$c_\alpha=\alpha(1-\alpha)/2$. Assume
\[
V^{\mathrm P}=\sup_{M\times(0,T]}v<\infty.
\]
Then, for $x\in M$ and $0<t\leq T$,
\begin{equation}
\label{eq:comparison-constant-curvature}
\alpha e^{-2Kt}\frac{|\nabla v|^2}{v}-\frac{v_t}{v}
\leq\frac{n(m-1)}{n(m-1)+\alpha(1+\alpha)e^{-2Kt}}\frac1t
+\frac{CK(m-1)V^{\mathrm P}}{\sqrt{c_\alpha}}.
\end{equation}

For the FDE, assume $1-2/n<m<1$ and fix $\alpha>1$.
Choose $\varepsilon$ and $\theta$ as in \eqref{eq:lambda-m-choice}, and set
\[
\delta_K=\min\left\{1,
\min_{\alpha\leq A\leq\alpha e^{2KT}}\lambda_{\min}(M(A))\right\}>0.
\]
Assume
\[
V^{\mathrm F}=\sup_{M\times(0,T]}(-v)<\infty.
\]
Then, for $x\in M$ and $0<t\leq T$,
\begin{equation}
\label{eq:comparison-constant-curvature-FDE}
\frac{v_t}{v}-\alpha e^{2K(T-t)}\frac{|\nabla v|^2}{v}
\leq\frac{C}{\delta_Kt}
+\frac{C\alpha K(1-m)e^{2KT}V^{\mathrm F}}{\delta_K}.
\end{equation}
Both estimates hold for every $K\geq0$.

Compared with the estimates in
\cite{LNVV09}, our gradient coefficients contain the additional time factor
$e^{-2Kt}$. 

\smallskip\noindent\textbf{Sharpness and coefficient choices.}
For the PME, $K=0$ gives
\[
\alpha\frac{|\nabla v|^2}{v}-\frac{v_t}{v}
\leq\frac{n(m-1)}{n(m-1)+\alpha(1+\alpha)}\frac1t.
\]
Letting $\alpha\uparrow1$ recovers the sharp bound $a_m/t$
in \cite{LNVV09}. 
For $K>0$, our PME bound and the estimate in \cite{HuangHuangLi2013}
have different forms and cannot be compared directly.

For the FDE, $K=0$  gives
\[
	\frac{v_t}{v}-\frac{|\nabla v|^2}{v}
	\leq\frac{C}{\delta_0t},
	\]
where $\delta_0$ is the value of $\delta_K$ at $K=0$.

Our normalized $L^p$ smallness condition applies to  collapsing manifolds  and allows localized negative curvature peaks, without
a volume lower bound or a uniform pointwise Ricci lower bound. Pointwise results
hold for every $K\geq0$ and may yield sharper constants or better
time-dependent coefficients.

	% ======================================================================

	\section{Geometric and analytic preliminaries}
	\label{sec:preliminaries}
	
	\subsection*{Geometric consequences of the curvature assumption}
	
	We collect the geometric consequences of the assumption
	$k(p,1)\leq\kappa$ that will be used below. Following Dai, Wei, and
	Zhang \cite{DWZ}, we use
	\[
	\rho_-(x)=\max\{0,-\lambda_{\min}(\Ric_x)\}
	\]
	to measure the negative part of the Ricci curvature. We first compare
	this convention with the Hilbert--Schmidt norm used by Zhang--Zhu
	\cite{ZhangZhu2017}. Let $\Ric^-$ be the negative spectral part of the
	Ricci endomorphism, and interpret its tensor norm as the
	Hilbert--Schmidt norm. Their normalized curvature quantity is
	\begin{equation}
		\label{eq:ZZ-curvature}
		k_{\mathrm{HS}}(p,r)
		=
		\sup_{x\in M}r^2
		\left(
		\fint_{B(x,r)}|\Ric^-|_{\mathrm{HS}}^p\,dV
		\right)^{1/p},
	\end{equation}
	and their smallness condition is
	\begin{equation}
		\label{1.1}
		k_{\mathrm{HS}}(p,1)\leq\kappa_{\mathrm{ZZ}}(n,p).
	\end{equation}
	
	To compare the two norms, let
	$\lambda_1(x)\leq\cdots\leq\lambda_n(x)$ be the Ricci eigenvalues and
	put $\lambda_i^-(x)=\max\{0,-\lambda_i(x)\}$. In an orthonormal
	eigenbasis,
	\[
	|\Ric_x^-|_{\mathrm{HS}}
	=\left(\sum_{i=1}^n(\lambda_i^-(x))^2\right)^{1/2},
	\qquad
	\rho_-(x)=\max_i\lambda_i^-(x).
	\]
	The pointwise comparison
	\eqref{eq:ricci-norm-comparison} implies
	\begin{equation}
		\label{eq:curvature-convention-comparison}
		k(p,r)\leq k_{\mathrm{HS}}(p,r)\leq\sqrt n\,k(p,r).
	\end{equation}
	The dimensional factor $\sqrt n$ can be absorbed into the initial
	smallness threshold $\kappa(n,p)$ by choosing
	$\kappa(n,p)\leq\kappa_{\mathrm{ZZ}}(n,p)/\sqrt n$.
	With this choice, our assumption implies \eqref{1.1}.
	
	We choose $\kappa=\kappa(n,p)>0$ sufficiently small so that the
	volume comparison, local Sobolev, cutoff, and heat kernel estimates
	in this section hold simultaneously.
	
	\subsection{Volume comparison, local Sobolev inequalities, and cutoff functions}

	\begin{lemma}[Volume and scale comparison]
		\label{lem:volume-scale}
		Let $n\geq3$ and $p>n/2$. If
		\[
		k(p,1)\leq\kappa,
		\]
		then, for every $x\in M$ and $0<r\leq R\leq1$,
		\begin{equation}
			\label{eq:volume-comparison}
			\frac{\vol(B(x,R))}{R^n}
			\leq
			C
			\frac{\vol(B(x,r))}{r^n}.
		\end{equation}
		Moreover,
		\begin{equation}
			\label{eq:scale-comparison}
			k(x,p,r)
			\leq
			C\,
			r^{2-n/p}k(x,p,1)
			\leq
			C\,
			r^{2-n/p}\kappa.
		\end{equation}
	\end{lemma}
	
	\begin{proof}
		Petersen--Wei volume comparison \cite{PW97} and local doubling \cite[Lemma~2.1]{ZhangZhu2017} give
		\eqref{eq:volume-comparison}.
		
		Since $B(x,r)\subset B(x,1)$,
		\begin{align*}
			k(x,p,r)
			&=
			r^2
			\left(
			\frac{1}{\vol(B(x,r))}
			\int_{B(x,r)}\rho_-^{\,p}\,dV
			\right)^{1/p}\\
			&\leq
			r^2
			\left(
			\frac{\vol(B(x,1))}{\vol(B(x,r))}
			\right)^{1/p}
			\left(
			\fint_{B(x,1)}\rho_-^{\,p}\,dV
			\right)^{1/p}.
		\end{align*}
		By \eqref{eq:volume-comparison} with $R=1$,
		\[
		\frac{\vol(B(x,1))}{\vol(B(x,r))}
		\leq Cr^{-n}.
		\]
		Consequently,
		\[
		k(x,p,r)
		\leq
		C
		r^{2-n/p}k(x,p,1).
		\]
		This proves \eqref{eq:scale-comparison}.
	\end{proof}
	
	The corresponding local Sobolev inequality is due to Dai, Wei, and Zhang.
	
	\begin{lemma}[Local Sobolev inequality]
		\label{lem:local-sobolev}
		Let $n\geq3$ and $p>n/2$. If
		\[
		k(p,1)\leq\kappa,
		\]
		then, for every $x\in M$, $0<r\leq1$, and
		$f\in C_0^\infty(B(x,r))$,
		\begin{equation}
			\label{eq:local-sobolev-dwz}
			\left(
			\fint_{B(x,r)}
			|f|^{\frac{2n}{n-2}}\,dV
			\right)^{\frac{n-2}{n}}
			\leq
			C\,r^2
			\fint_{B(x,r)}
			|\nabla f|^2\,dV.
		\end{equation}
	\end{lemma}
	
	\begin{proof}
		This is the normalized $L^2$-Sobolev inequality in
		\cite[Corollary~4.6]{DWZ}.
	\end{proof}

	The same curvature assumption also gives a cutoff with controlled Laplacian.
	
	\begin{corollary}[Cutoff function]
		\label{lem:cutoff}
		Let $n\geq3$ and $p>n/2$. Assume that
		\[
		k(p,1)\leq\kappa.
		\]
		Then, for every
		$x\in M$ and $0<r\leq1$, there exists a function
		$\phi\in C_0^\infty(B(x,r))$ satisfying
		\[
		0\leq\phi\leq1,
		\qquad
		\phi\equiv1\quad\text{on }B(x,r/2),
		\]
		and
		\begin{equation}
			\label{eq:powered-cutoff}
			|\Delta\phi|\leq Cr^{-2},
			\qquad
			|\nabla\phi|^2\leq Cr^{-2}\phi
			\quad\text{on }M,
		\end{equation}
	\end{corollary}
	
	\begin{proof}
		By \cite[Lemma~2.4]{ZhangZhu2017}, there exists
		$\psi\in C_0^\infty(B(x,r))$ such that
		\[
		0\leq\psi\leq1,
		\qquad
		\psi\equiv1\quad\text{on }B(x,r/2),
		\qquad
		|\nabla\psi|^2+|\Delta\psi|\leq Cr^{-2}.
		\]
		Set $\phi=\psi^2$, with the required support and values. Then
		\[
		|\nabla\phi|^2
		=4\psi^2|\nabla\psi|^2
		\leq Cr^{-2}\phi,
		\]
		and
		\[
		|\Delta\phi|
		=|2\psi\Delta\psi+2|\nabla\psi|^2|
		\leq Cr^{-2}.
		\]
		Thus $\phi$ satisfies \eqref{eq:powered-cutoff}.
	\end{proof}

	\subsection{Heat kernel estimates}
	\label{subsec:heat kernel}
	
	Throughout this paper, $\Delta$ denotes the Laplace--Beltrami operator
	\begin{equation*}
		\Delta=\operatorname{div}\nabla.
	\end{equation*}
	
	Let $G_\Omega(x,y,t)$ denote the Dirichlet heat kernel of $\Delta$
	on a smooth relatively compact domain $\Omega\Subset M$, and define
	\begin{equation*}
		P_t^\Omega f(x)
		=
		\int_\Omega
		G_\Omega(x,y,t)f(y)\,dV(y).
	\end{equation*}
	
	For a smooth domain $\Omega\Subset B(x_0,r)$, the Dirichlet
	Laplacian $\Delta_{\Omega,D}$ is self-adjoint and nonpositive.
	The Nash argument then gives the following heat kernel
	bound, with constants independent of $\Omega$ and $\vol(B(x_0,r))$.
	
	% ======================================================================
	\begin{lemma}[Dirichlet heat kernel upper bound]
		\label{lem:dirichlet-heat kernel}
		Let $n\geq3$ and $p>n/2$. Assume
		\[
		k(p,1)\leq\kappa.
		\]
		Then the following holds.
		
		Let $0<r\leq1$, $x_0\in M$, and let
		$\Omega\Subset B(x_0,r)$ be a smooth domain. Then, for every
		$t>0$,
		\begin{equation}
			\label{eq:dirichlet-kernel-upper}
			G_\Omega(z,y,t)
			\leq
			\frac{Cr^n}
			{\vol(B(x_0,r))\,t^{n/2}},
			\qquad z,y\in\Omega.
		\end{equation}
		If $p'=p/(p-1)$, then
		\begin{equation}
			\label{eq:dirichlet-kernel-Lp}
			\sup_{z\in\Omega}
			\|G_\Omega(z,\cdot,t)\|_{L^{p'}(\Omega)}
			\leq
			C
			\vol(B(x_0,r))^{-1/p}
			r^{n/p}t^{-n/(2p)}.
		\end{equation}
	\end{lemma}
	
	\begin{proof}
Set $\mathcal V_r=\vol(B(x_0,r))$.
Extend $f\in C_0^\infty(\Omega)$ by zero to $B(x_0,r)$.
The local Sobolev inequality gives
\begin{equation}\label{eq:unnormalized-local-sobolev}
\|f\|_{L^{2n/(n-2)}(\Omega)}^2
\leq C(r^n/\mathcal V_r)^{2/n}\|\nabla f\|_{L^2(\Omega)}^2.
\end{equation}
Interpolation between $L^1(\Omega)$ and $L^{2n/(n-2)}(\Omega)$ gives
\[
\|f\|_{L^2(\Omega)}^{2+4/n}
\leq\|f\|_{L^{2n/(n-2)}(\Omega)}^2\|f\|_{L^1(\Omega)}^{4/n}.
\]
Combining this with \eqref{eq:unnormalized-local-sobolev} yields
\begin{equation}\label{eq:local-Nash}
\|f\|_{L^2(\Omega)}^{2+4/n}
\leq C(r^n/\mathcal V_r)^{2/n}\|\nabla f\|_{L^2(\Omega)}^2\|f\|_{L^1(\Omega)}^{4/n}.
\end{equation}
By density, this holds on $H_0^1(\Omega)\cap L^1(\Omega)$.

For nonzero $f\in L^1(\Omega)\cap L^2(\Omega)$, set $u(t)=P_t^\Omega f$.
The energy identity and $L^1(\Omega)$ contraction give
\[
\begin{aligned}
\frac{d}{dt}\|u(t)\|_{L^2(\Omega)}^2
&=-2\|\nabla u(t)\|_{L^2(\Omega)}^2\\
&\leq-C^{-1}(\mathcal V_r/r^n)^{2/n}
\|f\|_{L^1(\Omega)}^{-4/n}\|u(t)\|_{L^2(\Omega)}^{2+4/n}.
\end{aligned}
\]
Integration yields
$\|u(t)\|_{L^2(\Omega)}^2\leq C(r^n/\mathcal V_r)t^{-n/2}\|f\|_{L^1(\Omega)}^2$.
Consequently,
\begin{equation}\label{eq:semigroup-L1-L2}
\|P_t^\Omega\|_{L^1(\Omega)\to L^2(\Omega)}
\leq C(r^n/\mathcal V_r)^{1/2}t^{-n/4}.
\end{equation}
For $f\in L^2(\Omega)$ and $h\in L^1(\Omega)\cap L^2(\Omega)$,
self-adjointness and \eqref{eq:semigroup-L1-L2} give
\[
\begin{aligned}
\left|\int_\Omega (P_t^\Omega f)h\,dV\right|
&=\left|\int_\Omega f(P_t^\Omega h)\,dV\right|\\
&\leq\|f\|_{L^2(\Omega)}\|P_t^\Omega h\|_{L^2(\Omega)}\\
&\leq C(r^n/\mathcal V_r)^{1/2}t^{-n/4}
\|f\|_{L^2(\Omega)}\|h\|_{L^1(\Omega)}.
\end{aligned}
\]
Taking suprema over $\|h\|_{L^1(\Omega)}\leq1$ and
$\|f\|_{L^2(\Omega)}\leq1$ yields
\[
\|P_t^\Omega\|_{L^2(\Omega)\to L^\infty(\Omega)}
\leq C(r^n/\mathcal V_r)^{1/2}t^{-n/4}.
\]
Applying these two bounds to
$P_t^\Omega=P_{t/2}^\Omega P_{t/2}^\Omega$ gives
\[
\|P_t^\Omega\|_{L^1(\Omega)\to L^\infty(\Omega)}
\leq C(r^n/\mathcal V_r)t^{-n/2}.
\]
The operator estimate gives \eqref{eq:dirichlet-kernel-upper}
for almost every $(z,y)\in\Omega\times\Omega$.
Since the heat kernel is continuous for $t>0$, the bound holds for every
$z,y\in\Omega$.

Since $\int_\Omega G_\Omega(z,y,t)\,dV(y)\leq1$, interpolation gives
\[
\|G_\Omega(z,\cdot,t)\|_{L^{p'}(\Omega)}
\leq\|G_\Omega(z,\cdot,t)\|_{L^\infty(\Omega)}^{1/p}
\leq C\mathcal V_r^{-1/p}r^{n/p}t^{-n/(2p)}.
\]
This proves \eqref{eq:dirichlet-kernel-Lp}, with constants independent of $\Omega$.
\end{proof}
	
	% ======================================================================
	% Standard remark for the volume factors
	% ======================================================================
	The following estimate requires no volume lower bound.

	% ======================================================================
	\begin{lemma}[Potential estimate]
	    \label{lem:kernel-potential}
	    Let $n\geq3$, $p>n/2$, and assume $k(p,1)\leq\kappa$.
	    Let $0<r\leq1$ and let $\Omega\Subset B(O,r)$ be a smooth domain.
	    For a nonnegative function $W\in L^p(\Omega)$, set
	    \[
	    \Lambda_W(r)=\frac{r^2\|W\|_{L^p(\Omega)}}{\vol(B(O,r))^{1/p}}.
	    \]
	    Then, for every $x\in\Omega$ and $t>0$,
	    \begin{equation}
	        \label{eq:potential-estimate}
	        P_t^\Omega W(x)
	        \leq C\Lambda_W(r)r^{-2}
	        \left(\frac{t}{r^2}\right)^{-\gamma},
	    \end{equation}
	    where $\gamma$ is defined in \eqref{eq:common-exponents}
	    and $C=C(n,p)>0$.
	\end{lemma}
	
	\begin{proof}
	    With $p'=p/(p-1)$, H\"older's inequality and
	    Lemma~\ref{lem:dirichlet-heat kernel} give
	    \begin{align*}
	        P_t^\Omega W(x)
	        &=\int_\Omega G_\Omega(x,y,t)W(y)\,dV(y)\\
	        &\leq\|G_\Omega(x,\cdot,t)\|_{L^{p'}(\Omega)}\|W\|_{L^p(\Omega)}\\
	        &\leq C\vol(B(O,r))^{-1/p}r^{n/p}t^{-n/(2p)}\|W\|_{L^p(\Omega)}\\
	        &=C\Lambda_W(r)r^{-2}\left(\frac{t}{r^2}\right)^{-\gamma}.
	    \end{align*}
	\end{proof}
	
	% ======================================================================

	\subsection{Negative Ricci curvature and smooth majorants}
	\label{subsec:negative-ricci}
	
	We shall apply Lemma~\ref{lem:kernel-potential} to a smooth majorant of
	$\rho_-$ supplied by the next lemma.

	\begin{lemma}[Smooth approximation of curvature]
		\label{lem:smooth-majorants}
		Let $\Omega\Subset M$ be a smooth domain.
		There exist an open neighborhood $U\supset\overline\Omega$ and
		functions $W_j\in C^\infty(U)$ such that
		\begin{equation*}
			W_j\geq\rho_-\geq0
			\quad\text{on }\overline\Omega,
			\qquad
			\|W_j-\rho_-\|_{L^\infty(\overline\Omega)}
			\longrightarrow0.
		\end{equation*}
		Consequently,
		\[
		W_j\to\rho_-
		\quad\text{in }L^q(\Omega)
		\quad\text{for every }1\leq q<\infty.
		\]
	\end{lemma}
	
	\begin{proof}
		Since $\mathrm{Ric}$ is smooth and the map
		\[
		A\longmapsto \max\{0,-\lambda_{\min}(A)\}
		\]
		is continuous on symmetric endomorphisms, $\rho_-$ is continuous.
		
		By compactness of $\overline\Omega$, smooth approximation gives $\widehat W_j\in C^\infty(U)$ with
		\[
		d_j
		=
		\|\widehat W_j-\rho_-\|_{L^\infty(\overline\Omega)}
		\longrightarrow0.
		\]
		Set
		\[
		W_j=\widehat W_j+d_j+j^{-1}.
		\]
		Then on $\overline\Omega$,
		\[
		W_j-\rho_-
		=
		(\widehat W_j-\rho_-)+d_j+j^{-1}
		\geq j^{-1}>0.
		\]
		Moreover,
		\[
		\|W_j-\rho_-\|_{L^\infty(\overline\Omega)}
		\leq 2d_j+j^{-1}
		\longrightarrow0.
		\]
		Uniform convergence and $\vol(\Omega)<\infty$ give convergence in $L^q(\Omega)$.
	\end{proof}
	
	% ======================================================================
	
	\subsection{Linearized auxiliary equation}
	\label{subsec:linearized}
	
	We next introduce the Dirichlet problem used to linearize the auxiliary
	equation.

	\begin{lemma}[Linearized auxiliary problem]
		\label{lem:duhamel-wt}
		Let $\Omega\Subset M$ be a smooth domain, let
		$W\in C^\infty(\overline\Omega)$ and assume $W\ge0$ on $\overline\Omega$.
		Let $T,\beta,b>0$. Let $h$ be the unique mild solution in
		$C([0,T];L^2(\Omega))$ of
		\begin{equation}
			\label{eq:h-problem}
			h_t=(\Delta_{\Omega,D}+2\beta W)h+2\beta bW,
			\qquad h(0)=0,
		\end{equation}
		and set $\zeta=h+b$. Then
		\begin{equation*}
			\zeta\in
			C(\overline\Omega\times[0,T])
			\cap C^\infty(\Omega\times(0,T])
		\end{equation*}
		satisfies
		\begin{equation}
			\label{eq:zeta-problem-general}
			\begin{cases}
				\zeta_t=\Delta \zeta+2\beta W\zeta,
				&\text{in }\Omega\times(0,T],\\
				\zeta(x,0)=b,
				&x\in\Omega,\\
				\zeta(x,t)=b,
				&(x,t)\in\partial\Omega\times(0,T].
			\end{cases}
		\end{equation}
		Moreover, on $\Omega\times(0,T]$, we have
		$\zeta\geq b$, $\zeta_t\geq0$, and
		\begin{equation}
			\label{eq:zeta-t-duhamel}
			\zeta_t(x,t)
			=
			2\beta bP_t^\Omega W(x)
			+
			2\beta\int_0^t\int_\Omega
			G_\Omega(x,y,t-s)W(y)\zeta_t(y,s)\,dV(y)\,ds.
		\end{equation}
	\end{lemma}
	
	\begin{proof}
		Since $W$ is bounded, the operator
		$A=\Delta_{\Omega,D}+2\beta W$ is self-adjoint with spectrum bounded above.
		By the spectral theorem, it generates an analytic semigroup
		$e^{tA}$ on $L^2(\Omega)$ \cite{Pazy1983}.
		The parabolic maximum principle gives $e^{tA}f\geq0$ for every
		nonnegative $f\in L^2(\Omega)$. The solution of \eqref{eq:h-problem}
		has the representation
		\[
		h(t)=\int_0^t e^{(t-s)A}(2\beta bW)\,ds.
		\]
		Thus $h\ge0$, $\zeta=h+b\ge b$, and standard Dirichlet parabolic
		regularity \cite{Lieberman1996} gives the asserted interior smoothness
		and boundary values. Since $W$, $\beta$, and $b$ do not depend on $t$,
		\[
		\zeta_t(x,t)=h_t(x,t)
		=\left[e^{tA}(2\beta bW)\right](x)\ge0.
		\]
		The function $\zeta_t$ has initial value $2\beta bW$ in
		$L^2(\Omega)$ and zero Dirichlet boundary values. Applying Duhamel's
		formula to its equation gives \eqref{eq:zeta-t-duhamel}.
	\end{proof}
	
	\subsection{Singular Volterra inequality}
	\label{subsec:volterra}
	
	The heat kernel estimates and the following Volterra inequality yield
	bounds for the auxiliary function depending on $W$ only through its
	normalized $L^p$ norm.
	
	% ======================================================================
	% Singular Volterra inequality
	% ======================================================================
	The following proposition is a weighted form of Henry's weakly singular
	Grönwall inequality \cite[Lemma~7.1.1]{Henry1981}. 
	
	\begin{proposition}[Singular Volterra inequality]
		\label{prop:volterra}
		Let $T>0$, $\gamma=n/(2p)\in(0,1)$, and $A,\lambda\geq0$. Define
		\begin{equation*}
			q_\gamma(t) = t^{-\gamma},
			\qquad t>0.
		\end{equation*}
		Let $f$ be a measurable, locally bounded map from $(0,T]$ to
		$[0,\infty)$.
		
		\medskip
		\noindent
		\textbf{(a) Singular forcing.}
		Assume
		\begin{equation*}
			\sup_{0<t\le T}\frac{f(t)}{q_\gamma(t)}<\infty
		\end{equation*}
		and
		\begin{equation}
			\label{eq:volterra-singular}
			f(t)
			\le
			Aq_\gamma(t)
			+
			\lambda\int_0^t q_\gamma(t-s)f(s)\,ds,
			\qquad 0<t\le T.
		\end{equation}
		Then, for some $C=C(n,p)>0$,
		\begin{equation}
			\label{eq:volterra-singular-conclusion}
			f(t)
			\le
			2Aq_\gamma(t)
			\exp\left\{
			(2C\lambda)^{1/(1-\gamma)}t
			\right\}.
		\end{equation}
		
		\medskip
		\noindent
		\textbf{(b) Constant forcing.}
		Assume instead that $f\in L^\infty(0,T]$ and
		\begin{equation}
			\label{eq:volterra-constant}
			f(t)
			\le
			A
			+
			\lambda\int_0^t q_\gamma(t-s)f(s)\,ds,
			\qquad 0<t\le T.
		\end{equation}
		Then
		\begin{equation}
			\label{eq:volterra-constant-conclusion}
			f(t)
			\le
			2A
			\exp\left\{
			(2C\lambda)^{1/(1-\gamma)}t
			\right\}.
		\end{equation}
	\end{proposition}
	
	\begin{proof}
The case $\lambda=0$ is immediate. Assume $\lambda>0$.
Let $\omega>0$.

\begin{equation}\label{eq:weighted-volterra-convolution}
\begin{aligned}
&\frac1{q_\gamma(t)}\int_0^t e^{-\omega(t-s)}q_\gamma(t-s)q_\gamma(s)\,ds\\
&\quad=\frac1{q_\gamma(t)}
\left(\int_0^{t/2}+\int_{t/2}^t\right)
e^{-\omega(t-s)}q_\gamma(t-s)q_\gamma(s)\,ds\\
&\quad\leq 2^\gamma e^{-\omega t/2}\int_0^{t/2}s^{-\gamma}\,ds
+2^\gamma\int_{t/2}^t e^{-\omega(t-s)}(t-s)^{-\gamma}\,ds\\
&\quad\leq\frac{2^{2\gamma-1}}{1-\gamma}e^{-\omega t/2}t^{1-\gamma}
+2^\gamma\int_0^\infty e^{-\omega u}u^{-\gamma}\,du\\
&\quad\leq C\omega^{\gamma-1}.
\end{aligned}
\end{equation}
The last inequality follows from
\[
\begin{aligned}
\int_0^\infty e^{-\omega u}u^{-\gamma}\,du
&=\omega^{\gamma-1}\int_0^\infty e^{-s}s^{-\gamma}\,ds
=\Gamma(1-\gamma)\omega^{\gamma-1},\\
\sup_{t>0}e^{-\omega t/2}t^{1-\gamma}
&=\left(\frac{2(1-\gamma)}{e}\right)^{1-\gamma}\omega^{\gamma-1}.
\end{aligned}
\]
Fix
\[
C=2^\gamma\Gamma(1-\gamma)
+\frac{2^{2\gamma-1}}{1-\gamma}
\left(\frac{2(1-\gamma)}{e}\right)^{1-\gamma},
\]
so that
\eqref{eq:weighted-volterra-convolution} holds, and choose
\[
\omega=(2C\lambda)^{1/(1-\gamma)}.
\]
Then
\begin{equation}\label{eq:volterra-key-inequality}
\lambda C\omega^{\gamma-1}=\tfrac12.
\end{equation}

For part~(a), fix $0<S\leq T$. By hypothesis, $M_S=\sup_{0<t\leq S}e^{-\omega t}f(t)/q_\gamma(t)<\infty$.
Multiply \eqref{eq:volterra-singular} by $e^{-\omega t}/q_\gamma(t)$.
Using \eqref{eq:weighted-volterra-convolution} gives
\[
e^{-\omega t}f(t)/q_\gamma(t)\leq A+M_S/2.
\]
Taking the supremum gives $M_S\leq2A$, proving part~(a).

For part~(b), set $N_S=\sup_{0<t\leq S}e^{-\omega t}f(t)<\infty$.
Multiply \eqref{eq:volterra-constant} by $e^{-\omega t}$ and use
\[
\lambda\int_0^\infty e^{-\omega u}u^{-\gamma}\,du
=\lambda\Gamma(1-\gamma)\omega^{\gamma-1}
\leq\lambda C\omega^{\gamma-1}=\tfrac12.
\]
It gives $N_S\leq A+N_S/2$, hence $N_S\leq2A$.
This proves part~(b).
\end{proof}
	
	% ======================================================================
	\section{The auxiliary function}
	\label{sec:auxiliary}
	
	We use $\Lambda(\beta)$, $G_r^{(\beta)}$, and
	$\underline J_r^{(\beta)}$ defined in
	\eqref{eq:Lambda-general}.
		To absorb the curvature term, we seek a positive function satisfying
	\[
	J_t-\Delta J+2WJ+a\frac{|\nabla J|^2}{J}=0.
	\]
	With $\beta=a-1$, the substitution
	$J=b^{1+1/\beta}\zeta^{-1/\beta}$ reduces this equation to
	$\zeta_t=\Delta\zeta+2\beta W\zeta$. Following
	\cite{ZhangZhu2017}, we estimate $\zeta$ by the Dirichlet heat kernel.
	We also require a bound for $\zeta_t/\zeta$ to control
	$\partial_t\log J$.

	\begin{proposition}[Auxiliary function construction]
		\label{prop:J}
		Assume $n\geq3$, $p>n/2$, and $k(p,1)\leq\kappa$.
		For $T>0$, $0<r\leq1$, $\beta,b>0$, and a smooth domain
		$\Omega\Subset B(O,r)$, there exist $W\in C^\infty(\overline\Omega)$
		and a positive function
		\[
		J\in C(\overline\Omega\times[0,T])
		\cap C^\infty(\Omega\times(0,T])
		\]
		such that
		\begin{equation}
			\label{eq:Lambda-W-bound}
			W\geq\rho_-\geq0\quad\text{on }\overline\Omega,
			\qquad
			\Lambda_W(r)=\frac{r^2\|W\|_{L^p(\Omega)}}{\vol(B(O,r))^{1/p}}
			\leq C\kappa.
		\end{equation}
		The function $J$ solves
		\begin{equation}
			\label{eq:auxiliary-equation-general}
			\begin{cases}
				J_t-\Delta J+2WJ
				+(\beta+1)\dfrac{|\nabla J|^2}{J}=0,
				&\text{in }\Omega\times(0,T],\\[6pt]
				J=b,
				&\text{on }(\overline\Omega\times\{0\})
				\cup(\partial\Omega\times(0,T]),
			\end{cases}
		\end{equation}
		and satisfies
		\begin{equation}
			\label{eq:J-bounds}
			b\,\underline J_r^{(\beta)}(t)\leq J(x,t)\leq b,
			\qquad
			-G_r^{(\beta)}(t)\leq\partial_t\log J(x,t)\leq0.
		\end{equation}
	\end{proposition}
	
\begin{proof}

By Lemma~\ref{lem:smooth-majorants}, choose $W_j\geq\rho_-$ converging
uniformly to $\rho_-$ on $\overline\Omega$. For sufficiently large $j$,
\[
r^2\vol(B(O,r))^{-1/p}\|W_j-\rho_-\|_{L^p(\Omega)}\leq\kappa.
\]
Fix $W=W_j$. Minkowski's inequality and Lemma~\ref{lem:volume-scale} give
\[
\Lambda_W(r)\leq k(O,p,r)+\kappa\leq C\kappa.
\]
Lemma~\ref{lem:duhamel-wt} supplies the solution of
\begin{equation}\label{eq:zeta-auxiliary-problem}
\zeta_t=\Delta\zeta+2\beta W\zeta,\qquad
\zeta=b\quad\text{on }(\overline\Omega\times\{0\})\cup
(\partial\Omega\times(0,T]).
\end{equation}
It satisfies $\zeta\geq b$ and $\zeta_t\geq0$. Define
\begin{equation}\label{eq:J-from-zeta}
J=b(\zeta/b)^{-1/\beta}.
\end{equation}
Differentiation gives \eqref{eq:auxiliary-equation-general} and
$\partial_t\log J=-\zeta_t/(\beta\zeta)\leq0$.

Set $f(t)=\sup_\Omega\zeta(\cdot,t)$. Duhamel's formula gives
\[
\zeta(x,t)=b+2\beta\int_0^t\int_\Omega
G_\Omega(x,y,t-s)W(y)\zeta(y,s)\,dV(y)\,ds.
\]
Taking the supremum over $x\in\Omega$ on both sides and applying
Lemma~\ref{lem:kernel-potential} gives
\begin{equation}\label{eq:f-volterra}
f(t)\leq b+C\beta\kappa r^{-2}
\int_0^t\left(\frac{t-s}{r^2}\right)^{-\gamma}f(s)\,ds.
\end{equation}
The bound $\zeta\leq b e^{2\beta\|W\|_\infty t}$ ensures that $f$ is bounded on finite intervals.
Apply Proposition~\ref{prop:volterra}(b) after rescaling time by $r^2$.
Since $1/(1-\gamma)=1+\sigma$, it gives
\[
f(t)\leq2b\exp\left(\beta\Lambda(\beta)t/r^2\right).
\]
By \eqref{eq:J-from-zeta},
\[
b\,2^{-1/\beta}e^{-\Lambda(\beta)t/r^2}\leq J\leq b.
\]

Set $h(t)=\sup_\Omega\zeta_t(\cdot,t)$. By
\eqref{eq:zeta-t-duhamel} and Lemma~\ref{lem:kernel-potential},
\begin{equation}\label{eq:H-volterra-singular}
\begin{aligned}
h(t)\leq{}&C\beta b\kappa r^{-2}(t/r^2)^{-\gamma}
+C\beta\kappa r^{-2}\int_0^t
((t-s)/r^2)^{-\gamma}h(s)\,ds.
\end{aligned}
\end{equation}
The identity $\zeta_t=e^{t(\Delta_{\Omega,D}+2\beta W)}(2\beta bW)$ gives
$h(t)\leq2\beta b\|W\|_\infty e^{2\beta\|W\|_\infty t}$.
Thus the weighted boundedness required by Proposition~\ref{prop:volterra}(a)
holds. Applying it after the same time rescaling gives
\begin{equation}\label{eq:zeta-t-bound}
0\leq\zeta_t(x,t)\leq C\beta b\kappa r^{-2}(t/r^2)^{-\gamma}
\exp\left(\beta\Lambda(\beta)t/r^2\right).
\end{equation}
Since $\zeta\geq b$ and $-\partial_t\log J=\zeta_t/(\beta\zeta)$,
we obtain the remaining bound in \eqref{eq:J-bounds}.
The constant $C=C(n,p)$ 
is chosen to cover both bounds.
\end{proof}

	% ======================================================================
\section{PME evolution inequalities and localization}
\label{sec:pme-evolution}

\subsection{Evolution inequalities}
\label{subsec:pme-evolution}
For $m>1$, define
\begin{equation}\label{eq:PME-operator}
\mathcal L=\partial_t-(m-1)v\Delta.
\end{equation}
The pressure equation and the two quantities used below are
\begin{equation}\label{eq:z-q-PME}
v_t=(m-1)v\Delta v+|\nabla v|^2,\qquad
z=\frac{v_t}{v},\qquad q_{\mathrm P}=\frac{|\nabla v|^2}{v}\geq0.
\end{equation}
In particular,
\begin{equation}\label{eq:pressure-identities}
z-q_{\mathrm P}=(m-1)\Delta v.
\end{equation}
Differentiation and the Bochner formula give
\begin{align}
\mathcal Lz&=2m\langle\nabla z,\nabla v\rangle+z^2,
\label{eq:Lz}\\
\mathcal Lq_{\mathrm P}
&=2m\langle\nabla q_{\mathrm P},\nabla v\rangle
+2zq_{\mathrm P}-q_{\mathrm P}^2
-2(m-1)\bigl(|\nabla^2v|^2+\Ric(\nabla v,\nabla v)\bigr).
\label{eq:Lq}
\end{align}
These are the identities in \cite{LNVV09}.
For a positive smooth function $J$, define
\begin{equation}\label{eq:F-PME-definition}
F_{\mathrm P}=\alpha Jq_{\mathrm P}-z.
\end{equation}

\begin{lemma}[Evolution inequality for $F_{\mathrm P}$]
\label{lem:evolution-F}
Let $0<\alpha<1$, let $J$ be positive and smooth, and let
$V_0\geq0$ be smooth with $\Ric\geq-V_0g$. Then
\begin{equation}\label{eq:evolution-F-pme}
\begin{aligned}
\mathcal LF_{\mathrm P}\leq{}&2m\langle\nabla F_{\mathrm P},\nabla v\rangle
-z^2+2\alpha Jzq_{\mathrm P}-\frac{\alpha(1+\alpha)}2Jq_{\mathrm P}^2
-\alpha(1+\alpha)(m-1)J|\nabla^2v|^2\\
&+\alpha(m-1)\left[J_t-\Delta J+2V_0J+
 a_\alpha^{\mathrm{PME}}\frac{|\nabla J|^2}{J}\right]|\nabla v|^2
+\alpha J_t(1-(m-1)v)q_{\mathrm P}.
\end{aligned}
\end{equation}
\end{lemma}
\begin{proof}
Using \eqref{eq:Lz}--\eqref{eq:Lq} and
$v\nabla q_{\mathrm P}=\nabla(|\nabla v|^2)-q_{\mathrm P}\nabla v$, we obtain
\begin{equation}\label{eq:LF-raw}
\begin{aligned}
\mathcal LF_{\mathrm P}={}&2m\langle\nabla F_{\mathrm P},\nabla v\rangle
-z^2+2\alpha Jzq_{\mathrm P}-\alpha Jq_{\mathrm P}^2\\
&+\alpha(J_t-(m-1)v\Delta J)q_{\mathrm P}
-2\alpha(m-1)J\bigl(|\nabla^2v|^2+\Ric(\nabla v,\nabla v)\bigr)\\
&-4\alpha(m-1)\nabla^2v(\nabla J,\nabla v)
-2\alpha q_{\mathrm P}\langle\nabla J,\nabla v\rangle.
\end{aligned}
\end{equation}
Young's inequality gives
\begin{align*}
-4\alpha(m-1)\nabla^2v(\nabla J,\nabla v)
&\leq\alpha(m-1)(1-\alpha)J|\nabla^2v|^2
+\frac{4\alpha(m-1)}{1-\alpha}\frac{|\nabla J|^2}{J}|\nabla v|^2,\\
-2\alpha q_{\mathrm P}\langle\nabla J,\nabla v\rangle
&\leq c_\alpha Jq_{\mathrm P}^2
+\frac{2\alpha}{1-\alpha}\frac{|\nabla J|^2}{J}|\nabla v|^2.
\end{align*}
The coefficients of $\frac{|\nabla J|^2}{J}|\nabla v|^2$ sum to $\alpha(m-1)a_\alpha^{\mathrm{PME}}$.
The Ricci lower bound and the identity
\[
\begin{aligned}
\alpha(J_t-(m-1)v\Delta J)q_{\mathrm P}
={}&\alpha(m-1)(J_t-\Delta J)|\nabla v|^2+\alpha J_t(1-(m-1)v)q_{\mathrm P}
\end{aligned}
\]
give \eqref{eq:evolution-F-pme}.
\end{proof}

Take $W,J$ from Proposition~\ref{prop:J} with
$\beta=\beta^{\mathrm{PME}}$ and $b=1$.
Set $V_0=W$ and $G=G_{r,\mathrm P}$. Then
$0<J\leq1$ and $-G(t)\leq\partial_t\log J\leq0$.
In \eqref{eq:evolution-F-pme}, the auxiliary equation gives
\begin{equation}\label{eq:auxiliary-J}
\alpha(m-1)\left[J_t-\Delta J+2V_0J+
 a_\alpha^{\mathrm{PME}}\frac{|\nabla J|^2}{J}\right]|\nabla v|^2=0.
\end{equation}
Since $J_t\leq0$ and $-J_t\leq GJ$,
\[
\alpha J_t(1-(m-1)v)q_{\mathrm P}
\leq\alpha(m-1)G(t)vJq_{\mathrm P}.
\]
Since $z=\alpha Jq_{\mathrm P}-F_{\mathrm P}$ and $0<J\leq1$,
\[
\begin{aligned}
&-z^2+2\alpha Jzq_{\mathrm P}-\frac{\alpha(1+\alpha)}2Jq_{\mathrm P}^2\\
&\quad=-F_{\mathrm P}^2-c_\alpha Jq_{\mathrm P}^2
-\alpha^2J(1-J)q_{\mathrm P}^2
\leq-F_{\mathrm P}^2-c_\alpha Jq_{\mathrm P}^2.
\end{aligned}
\]
The trace inequality and \eqref{eq:pressure-identities} give
\begin{equation}\label{eq:PME-Hessian-coercivity}
-\alpha(1+\alpha)(m-1)J|\nabla^2v|^2
\leq-\frac{\alpha(1+\alpha)J}{n(m-1)}(z-q_{\mathrm P})^2.
\end{equation}
Substituting these estimates into \eqref{eq:evolution-F-pme} yields
\begin{equation}\label{eq:F-evolution-PME}
\begin{aligned}
\mathcal LF_{\mathrm P}\leq{}&2m\langle\nabla F_{\mathrm P},\nabla v\rangle
-F_{\mathrm P}^2-c_\alpha Jq_{\mathrm P}^2
-\frac{\alpha(1+\alpha)J}{n(m-1)}(z-q_{\mathrm P})^2
+\alpha(m-1)G(t)vJq_{\mathrm P}.
\end{aligned}
\end{equation}

\subsection{Localized maximum principle}
\label{subsec:pme-localization}
\begin{lemma}[Localized maximum principle for the PME]
\label{lem:localized-maximum}
Let $\Omega\Subset B(O,r)$ contain $\overline{B(O,r/2)}$, and choose
a cutoff function $\phi\in C_0^\infty(\Omega)$ as in
Corollary~\ref{lem:cutoff}, so that
\begin{equation}\label{eq:cutoff-properties}
0\leq\phi\leq1,\qquad\phi=1\text{ on }B(O,r/2),
\end{equation}
\begin{equation}\label{eq:cutoff-bounds}
|\Delta\phi|\leq Cr^{-2},\qquad |\nabla\phi|^2\leq Cr^{-2}\phi.
\end{equation}
Suppose \eqref{eq:F-evolution-PME} holds, $0<J_0(t)\leq J\leq1$,
$J_0$ is nonincreasing, $G\geq0$, and $tG(t)$ is nondecreasing.
Then on $B(O,r/2)\times(0,T]$,
\begin{equation}\label{eq:localized-maximum-conclusion-PME}
\begin{aligned}
F_{\mathrm P}\leq{}&
\frac{n(m-1)}{n(m-1)+\alpha(1+\alpha)J_0(t)}\frac1t\\
&+\frac{Cm^2(m-1)V_r^{\mathrm P}(t)}{c_\alpha J_0(t)r^2}
+\frac{C(m-1)G(t)V_r^{\mathrm P}(t)}{\sqrt{c_\alpha}}.
\end{aligned}
\end{equation}
\end{lemma}
\begin{proof}
Fix $(x_*,t_*)$ in the smaller cylinder and $0<\tau<t_*$.
If the maximum of $(s-\tau)\phi^2F_{\mathrm P}$ on $\overline\Omega\times[\tau,t_*]$ is nonpositive,
the conclusion is immediate. Otherwise, let $(x_0,s_0)$ be a point where this maximum is attained.
Let $X=\phi F_{\mathrm P}$,
$R=\phi q_{\mathrm P}$, and $V=V_r^{\mathrm P}(s_0)$.
Then $\nabla F_{\mathrm P}=-2F_{\mathrm P}\nabla\phi/\phi$.
Using \eqref{eq:F-evolution-PME}, \eqref{eq:cutoff-bounds}, $s_0>\tau$,
$\phi\leq1$, and $|\nabla v|^2\leq Vq_{\mathrm P}$, we obtain
\[
\begin{aligned}
0\leq{}&\frac{X}{s_0-\tau}-X^2-c_\alpha JR^2
-\frac{\alpha(1+\alpha)J}{n(m-1)}[X+(1-\alpha J)R]^2\\
&+\frac{Cm\sqrt V}{r}X\sqrt R
+\frac{C(m-1)V}{r^2}X
+\alpha(m-1)G(s_0)VJ\phi R.
\end{aligned}
\]

Since $1-\alpha J\geq1-\alpha>0$,
\[
\frac{Cm\sqrt V}{r}X\sqrt R
\leq \frac{2\alpha(1+\alpha)J(1-\alpha J)}{n(m-1)}XR
+\frac{Cm^2(m-1)V}{c_\alpha J_0(s_0)r^2}X.
\]
Young's inequality also gives
\[
\alpha(m-1)G(s_0)VJ\phi R
\leq\frac{c_\alpha J}{2}R^2
+\frac{C(m-1)^2G(s_0)^2V^2}{c_\alpha}.
\]
Combining the terms linear in $X$, expanding the square, and using $J\geq J_0(s_0)$ gives
\[
\left(1+\frac{\alpha(1+\alpha)J_0(s_0)}{n(m-1)}\right)X^2
\leq\left(\frac1{s_0-\tau}+
\frac{Cm^2(m-1)V}{c_\alpha J_0(s_0)r^2}\right)X
+\frac{C(m-1)^2G(s_0)^2V^2}{c_\alpha}.
\]
Solving this quadratic inequality gives
\[
\begin{aligned}
X\leq{}&\frac{n(m-1)}{n(m-1)+\alpha(1+\alpha)J_0(s_0)}\frac1{s_0-\tau}\\
&+\frac{Cm^2(m-1)V}{c_\alpha J_0(s_0)r^2}
+\frac{C(m-1)G(s_0)V}{\sqrt{c_\alpha}}.
\end{aligned}
\]

Maximality gives $(t_*-\tau)F_{\mathrm P}(x_*,t_*)\leq(s_0-\tau)X$.
Use $s_0-\tau\leq s_0\leq t_*$, $V\leq V_r^{\mathrm P}(t_*)$,
$J_0(s_0)\geq J_0(t_*)$, and $s_0G(s_0)\leq t_*G(t_*)$.
Divide by $t_*-\tau$ and let $\tau\downarrow0$.
\end{proof}

We now deduce Theorem~\ref{thm:local-AB} from Lemma~\ref{lem:localized-maximum}.
\begin{proof}[Proof of \cref{thm:local-AB}]
Compactness of $\overline{B(O,r)}\times[0,t]$ gives $V_r^{\mathrm P}(t)<\infty$.
Choose $\phi$ as in Corollary~\ref{lem:cutoff} and a smooth domain with
$\overline{B(O,r/2)}\cup\supp\phi\Subset\Omega\Subset B(O,r)$.
Proposition~\ref{prop:J}, with $b=1$ and $\beta=\beta^{\mathrm{PME}}$, gives
\[
\underline J_{r,\mathrm P}\leq J\leq1,\qquad
-G_{r,\mathrm P}\leq\partial_t\log J\leq0.
\]
Apply \eqref{eq:F-evolution-PME} and Lemma~\ref{lem:localized-maximum}
with $V_0=W$, $J_0=\underline J_{r,\mathrm P}$, and $G=G_{r,\mathrm P}$.
Since $tG_{r,\mathrm P}(t)$ is nondecreasing and
$\alpha\underline J_{r,\mathrm P}q_{\mathrm P}-z\leq F_{\mathrm P}$,
we obtain \eqref{eq:main-estimate-PME}.
\end{proof}

\section{FDE evolution inequalities and localization}
\label{sec:fde-evolution}

\subsection{Evolution inequalities}
\label{subsec:fde-evolution}
For $1-2/n<m<1$, define
\begin{equation}\label{eq:pressure-FDE-evolution}
v=\frac{m}{m-1}u^{m-1}<0,\qquad
\mathcal L=\partial_t-(m-1)v\Delta.
\end{equation}
Set
\begin{equation}\label{eq:fde-z-q}
z=\frac{v_t}{v},\qquad q_{\mathrm F}=-\frac{|\nabla v|^2}{v}\geq0.
\end{equation}
The pressure equation implies
\begin{equation}\label{eq:fde-pressure-identity}
z+q_{\mathrm F}=(m-1)\Delta v,\qquad
|\nabla v|^2=(-v)q_{\mathrm F}.
\end{equation}
As in the PME calculation,
\begin{align}
\mathcal Lz&=2m\langle\nabla z,\nabla v\rangle+z^2,
\label{eq:Lz-fde}\\
\mathcal Lq_{\mathrm F}
&=2m\langle\nabla q_{\mathrm F},\nabla v\rangle
+2zq_{\mathrm F}+q_{\mathrm F}^2
+2(m-1)\bigl(|\nabla^2v|^2+\Ric(\nabla v,\nabla v)\bigr).
\label{eq:Lq-fde}
\end{align}
For a positive smooth function $J$, define directly
\begin{equation}\label{eq:fde-F-explicit}
F_{\mathrm F}=z+Jq_{\mathrm F}
=\frac{v_t}{v}-J\frac{|\nabla v|^2}{v}.
\end{equation}

\begin{lemma}[Evolution inequality for $F_{\mathrm F}$]
\label{lem:evolution-F-fde}
Let $J$ be positive and smooth, and let $V_0\geq0$ be smooth with
$\Ric\geq-V_0g$. Use $\mu$, $\varepsilon$, $\theta$, and
$a^{\mathrm{FDE}}$ from \eqref{eq:lambda-m-choice}--\eqref{eq:a-FDE}.
Then
\begin{equation}\label{eq:evolution-F-fde}
\begin{aligned}
\mathcal LF_{\mathrm F}\leq{}&2m\langle\nabla F_{\mathrm F},\nabla v\rangle
+z^2+2Jzq_{\mathrm F}+(1+\theta)Jq_{\mathrm F}^2
-2\mu(1-\varepsilon)J|\nabla^2v|^2\\
&+\mu\left[J_t-\Delta J+2V_0J+
a^{\mathrm{FDE}}\frac{|\nabla J|^2}{J}\right]|\nabla v|^2
+J_t(1-(m-1)v)q_{\mathrm F}.
\end{aligned}
\end{equation}
\end{lemma}
\begin{proof}
From \eqref{eq:Lz-fde}--\eqref{eq:Lq-fde} and
$v\nabla q_{\mathrm F}=-\nabla(|\nabla v|^2)-q_{\mathrm F}\nabla v$,
\[
\begin{aligned}
\mathcal LF_{\mathrm F}={}&2m\langle\nabla F_{\mathrm F},\nabla v\rangle
+z^2+2Jzq_{\mathrm F}+Jq_{\mathrm F}^2+(J_t-(m-1)v\Delta J)q_{\mathrm F}\\
&-2\mu J\bigl(|\nabla^2v|^2+\Ric(\nabla v,\nabla v)\bigr)
-4\mu\nabla^2v(\nabla J,\nabla v)
-2q_{\mathrm F}\langle\nabla J,\nabla v\rangle.
\end{aligned}
\]

Young's inequality and the Ricci lower bound give
\begin{align*}
-4\mu\nabla^2v(\nabla J,\nabla v)
&\leq2\mu\varepsilon J|\nabla^2v|^2
+\frac{2\mu}{\varepsilon}\frac{|\nabla J|^2}{J}|\nabla v|^2,\\
-2q_{\mathrm F}\langle\nabla J,\nabla v\rangle
&\leq\theta Jq_{\mathrm F}^2
+\frac1\theta\frac{|\nabla J|^2}{J}|\nabla v|^2,\\
-2\mu J\Ric(\nabla v,\nabla v)&\leq2\mu V_0J|\nabla v|^2.
\end{align*}
The coefficients of $\frac{|\nabla J|^2}{J}|\nabla v|^2$ sum to $\mu a^{\mathrm{FDE}}$.
Also,
\[
(J_t-(m-1)v\Delta J)q_{\mathrm F}
=\mu(J_t-\Delta J)|\nabla v|^2
+J_t(1-(m-1)v)q_{\mathrm F}.
\]
Combining these identities proves \eqref{eq:evolution-F-fde}.
\end{proof}

Take $W,J$ from Proposition~\ref{prop:J} with
$\beta=\beta^{\mathrm{FDE}}$ and boundary value $b=\alpha/\underline J_{r,\mathrm F}(T)$.
Set $V_0=W$ and $G=G_{r,\mathrm F}$. Then
$\alpha\leq J\leq b$ and $-G(t)\leq\partial_t\log J\leq0$.
In \eqref{eq:evolution-F-fde}, the auxiliary equation gives
\begin{equation}\label{eq:fde-auxiliary-equation}
\mu\left[J_t-\Delta J+2V_0J+
a^{\mathrm{FDE}}\frac{|\nabla J|^2}{J}\right]|\nabla v|^2=0.
\end{equation}
Using $J_t\leq0$ and $-J_t\leq GJ$, we obtain
\begin{equation}\label{eq:fde-Jt-bound}
J_t(1-(m-1)v)q_{\mathrm F}
\leq(1-m)G(t)(-v)Jq_{\mathrm F}.
\end{equation}
The trace inequality and \eqref{eq:fde-pressure-identity} give
\begin{equation}\label{eq:fde-hessian-coercivity}
-2\mu(1-\varepsilon)J|\nabla^2v|^2
\leq-\lambda_m J(z+q_{\mathrm F})^2.
\end{equation}
For $A=J$ and $z=F_{\mathrm F}-Aq_{\mathrm F}$, the remaining quadratic expression is
\begin{equation}\label{eq:fde-matrix}
z^2+2Azq_{\mathrm F}+A(1+\theta)q_{\mathrm F}^2
-\lambda_m A(z+q_{\mathrm F})^2
=-\begin{pmatrix}F_{\mathrm F}&q_{\mathrm F}\end{pmatrix}
M(A)\begin{pmatrix}F_{\mathrm F}\\q_{\mathrm F}\end{pmatrix}.
\end{equation}

\Needspace{10\baselineskip}
\begin{lemma}[Positive definiteness of $M(A)$]
\label{lem:matrix-M-positive}
Let $A_0>1$ and $0<\theta<1$ satisfy
\[
A_0(1-\theta)-1>0,\qquad
\lambda_m>\frac{A_0-1-\theta}{A_0(1-\theta)-1}.
\]
Then the symmetric matrix
\begin{equation}\label{eq:fde-coercive-matrix}
M(A)=\begin{pmatrix}
\lambda_m A-1&\lambda_m A(1-A)\\
\lambda_m A(1-A)&A[\lambda_m(1-A)^2-(1+\theta-A)]
\end{pmatrix}
\end{equation}
is positive definite for every $A\geq A_0$.
\end{lemma}
\begin{proof}
For $A\geq A_0>1$, the assumptions give
\[
\begin{aligned}
&\lambda_m
>\frac{A_0-1-\theta}{A_0(1-\theta)-1}
=\frac1{1-\theta}
+\frac{\theta^2}{(1-\theta)[A_0(1-\theta)-1]}
>\frac1{1-\theta},\\
&M_{11}(A)
=\lambda_m A-1>\frac{A}{1-\theta}-1>0.
\end{aligned}
\]
For $A\geq A_0$, the assumptions also give
\[
\begin{aligned}
\det M(A)
&=A\{A[\lambda_m(1-\theta)-1]+1+\theta-\lambda_m\}\\
&\geq A\{A_0[\lambda_m(1-\theta)-1]+1+\theta-\lambda_m\}>0.
\end{aligned}
\]
Since $M_{11}(A)>0$ and $\det M(A)>0$, the symmetric matrix $M(A)$
is positive definite for every $A\geq A_0$.
\end{proof}
The choice \eqref{eq:lambda-m-choice} implies
\[
\alpha(1-\theta)-1>0,\qquad
\lambda_m[\alpha(1-\theta)-1]>\alpha-1-\theta.
\]
Apply Lemma~\ref{lem:matrix-M-positive} with $A_0=\alpha$.
By \eqref{eq:delta-effective-theorem}, $M(A)\geq\delta I$ for
$\alpha\leq A\leq b$. Equations \eqref{eq:evolution-F-fde}--\eqref{eq:fde-matrix}
therefore give
\begin{equation}\label{eq:F-evolution-fde-coercive}
\mathcal LF_{\mathrm F}\leq2m\langle\nabla F_{\mathrm F},\nabla v\rangle
-\delta(F_{\mathrm F}^2+q_{\mathrm F}^2)
+b(1-m)G(t)(-v)q_{\mathrm F}.
\end{equation}

\subsection{Localized maximum principle}
\label{subsec:fde-localization}
\begin{lemma}[Localized maximum principle for the FDE]
\label{lem:fde-localized}
Let $\Omega\Subset B(O,r)$ contain $\overline{B(O,r/2)}$, and let
$\phi\in C_0^\infty(\Omega)$ satisfy
\eqref{eq:cutoff-properties}--\eqref{eq:cutoff-bounds}.
Suppose \eqref{eq:F-evolution-fde-coercive} holds, $0<\delta\leq1$,
$G\geq0$, and $tG(t)$ is nondecreasing. Then on $B(O,r/2)\times(0,T]$,
\begin{equation}\label{eq:fde-localized}
F_{\mathrm F}\leq\frac{C}{\delta t}+\frac{Cm^2V_r^{\mathrm F}(t)}{\delta^2r^2}+\frac{C(1-m)V_r^{\mathrm F}(t)}{\delta r^2}+\frac{Cb(1-m)G(t)V_r^{\mathrm F}(t)}{\delta}.
\end{equation}
\end{lemma}
\begin{proof}
Fix $(x_*,t_*)$ and $0<\tau<t_*$.
At a positive maximum $(x_0,s_0)$ of $(s-\tau)\phi^2F_{\mathrm F}$,
let $X=\phi F_{\mathrm F}$, $R=\phi q_{\mathrm F}$,
and $V=V_r^{\mathrm F}(s_0)$. Using \eqref{eq:F-evolution-fde-coercive},
\eqref{eq:cutoff-bounds}, $|\nabla v|^2\leq Vq_{\mathrm F}$, and $\phi\leq1$, we obtain
\begin{equation}\label{eq:fde-local-maximum-raw}
0\leq\frac{X}{s_0-\tau}-\delta(X^2+R^2)+\frac{Cm\sqrt V}{r}X\sqrt R+\frac{C(1-m)V}{r^2}X+b(1-m)G(s_0)VR.
\end{equation}

Young's inequality gives
\begin{align*}
\frac{Cm\sqrt V}{r}X\sqrt R
&\leq\frac{\delta}{4}(X^2+R^2)
+\frac{Cm^4V^2}{\delta^3r^4},\\
\frac{C(1-m)V}{r^2}X
&\leq\frac{\delta}{4}X^2
+\frac{C(1-m)^2V^2}{\delta r^4},\\
b(1-m)G(s_0)VR
&\leq\frac{\delta}{4}R^2
+\frac{Cb^2(1-m)^2G(s_0)^2V^2}{\delta}.
\end{align*}
Substitution into \eqref{eq:fde-local-maximum-raw} and the quadratic formula yield
\[
X\leq\frac{C}{(s_0-\tau)\delta}
+\frac{Cm^2V}{\delta^2r^2}
+\frac{C(1-m)V}{\delta r^2}
+\frac{Cb(1-m)G(s_0)V}{\delta}.
\]

Maximality gives $(t_*-\tau)F_{\mathrm F}(x_*,t_*)\leq(s_0-\tau)X$.
Use $s_0-\tau\leq s_0\leq t_*$, $V\leq V_r^{\mathrm F}(t_*)$,
and $s_0G(s_0)\leq t_*G(t_*)$. Divide by $t_*-\tau$
and let $\tau\downarrow0$ to obtain \eqref{eq:fde-localized}.
\end{proof}

We now deduce Theorem~\ref{thm:fde-local-AB} from Lemma~\ref{lem:fde-localized}.
\begin{proof}[Proof of \cref{thm:fde-local-AB}]
Positivity and continuity of $u$ on the compact cylinder
$\overline{B(O,r)}\times[0,t]$ give $V_r^{\mathrm F}(t)<\infty$.
Choose $\phi$ and $\Omega$ as in the PME proof.
Proposition~\ref{prop:J}, with $\beta=\beta^{\mathrm{FDE}}$ and boundary value $b$, gives
\begin{equation}\label{eq:J-bounds-FDE-proof}
b\underline J_{r,\mathrm F}(t)\leq J\leq b,\qquad
-G_{r,\mathrm F}(t)\leq\partial_t\log J\leq0.
\end{equation}
Since $b\underline J_{r,\mathrm F}(t)\geq b\underline J_{r,\mathrm F}(T)=\alpha$,
Lemmas~\ref{lem:evolution-F-fde} and \ref{lem:matrix-M-positive} yield
\eqref{eq:F-evolution-fde-coercive} with $V_0=W$.
Lemma~\ref{lem:fde-localized} applies with $G=G_{r,\mathrm F}$,
as $tG_{r,\mathrm F}(t)$ is nondecreasing.
Together with
\[
z+b\underline J_{r,\mathrm F}(t)q_{\mathrm F}\leq F_{\mathrm F},
\]
this proves \eqref{eq:fde-theorem-estimate}.
\end{proof}

\section{Integrated Harnack inequalities}
\label{sec:unified-harnack}
In both proofs, set $\xi(s)=(1-s)t_1+st_2$ along the curve
$\eta:[0,1]\to B(O,r/2)$ in the corresponding theorem.

\begin{proof}[Proof of Theorem~\ref{thm:integrated-Harnack}]
By \eqref{eq:main-estimate-PME},
\begin{equation}\label{eq:AB-rearranged-PME}
\frac{v_t}{v}\geq\alpha\underline J_{r,\mathrm P}(t)
\frac{|\nabla v|^2}{v}
-\frac{n(m-1)}{n(m-1)+\alpha(1+\alpha)\underline J_{r,\mathrm P}(t)}\frac1t
-\Psi_r(t).
\end{equation}
Along $(\eta(s),\xi(s))$, completing the square gives
\[
\begin{aligned}
\frac{d}{ds}\log v
&=\langle\nabla\log v,\dot\eta\rangle+(t_2-t_1)\,\frac{v_t}{v}\\
&\geq-\frac{|\dot\eta|^2}
{4\alpha(t_2-t_1)\,\underline J_{r,\mathrm P}(\xi)v}
-(t_2-t_1)\left[
\frac{n(m-1)}{n(m-1)+\alpha(1+\alpha)\underline J_{r,\mathrm P}(\xi)}\frac1{\xi}
+\Psi_r(\xi)\right].
\end{aligned}
\]

Integrating over $[0,1]$ and exponentiating gives \eqref{eq:harnack-PME}.
For Corollary~\ref{cor:distance-harnack}, continuity and positivity of $v$ on the compact cylinder imply $v_{*,\rho}>0$.
The corollary follows by integrating along curves in $B(O,\rho)$
and taking the infimum over these curves, using $v\geq v_{*,\rho}$
and the monotonicity of $\underline J_{r,\mathrm P}$.
\end{proof}

\begin{proof}[Proof of Theorem~\ref{thm:fde-harnack}]
Set $w=-v>0$. Equation \eqref{eq:fde-theorem-estimate} gives
\begin{equation}\label{eq:AB-rearranged-FDE}
\frac{w_t}{w}\leq\Phi_r(t)
-\alpha\frac{\underline J_{r,\mathrm F}(t)}{\underline J_{r,\mathrm F}(T)}w|\nabla\log w|^2.
\end{equation}
Along $(\eta(s),\xi(s))$, completion of the square yields
\[
\begin{aligned}
\frac{d}{ds}\log w
&=\langle\nabla\log w,\dot\eta\rangle+(t_2-t_1)\,\frac{w_t}{w}\\
&\leq\frac{\underline J_{r,\mathrm F}(T)|\dot\eta|^2}
{4\alpha(t_2-t_1)\,\underline J_{r,\mathrm F}(\xi)w}
+(t_2-t_1)\,\Phi_r(\xi).
\end{aligned}
\]
Integrating over $[0,1]$ and exponentiating proves \eqref{eq:fde-harnack}.
For Corollary~\ref{cor:fde-distance-harnack}, continuity and positivity of $-v$ on the compact cylinder imply $v_{*,\rho}^{\mathrm F}>0$.
The corollary follows by integrating along curves in $B(O,\rho)$
and taking the infimum over these curves, using $-v\geq v_{*,\rho}^{\mathrm F}$
and the monotonicity of $\underline J_{r,\mathrm F}$.
\end{proof}

\end{document}